\def\l{\lambda}

\def\qed{$\Box$\medskip}

\newtheorem{theorem}{Theorem} \newtheorem{lemma}[theorem]{Lemma}

\newtheorem{corollary}[theorem]{Corollary}
\newtheorem{definition}[theorem]{Definition}

\def\proof{\noindent{\bf Proof.}\hspace{2mm}}
\documentclass[12pt]{article} %\usepackage{times} %\usepackage{mathart}
\usepackage{amsfonts} \usepackage{latexsym} \usepackage{amssymb,amsmath}

\author{Juliette Kennedy \thanks{Research partially
supported by grant 40734 of the Academy of Finland.}\\ Department of
Mathematics\\ University of Helsinki\\ Helsinki, Finland \and Saharon
Shelah \thanks{The second author would like to thank the Israel
Science Foundation for
      partial support of this research (Grant no. 242/03). Publication 852.}\\
Institute of Mathematics\\ Hebrew University\\ Jerusalem,
Israel\\} \title{More on regular reduced products}

\begin{document}

\maketitle

\begin{abstract}
The authors show, by means of a finitary version $\square^{fin}_{\lambda, D}$ of the combinatorial
principle $\square^{b^*}_{\lambda}$ of \cite{269}, the consistency of the failure, relative to the
consistency of supercompact cardinals, of the following: for all regular filters $D$ on a cardinal
$\lambda$, if $M_i$ and $N_i$ are elementarily equivalent models of a language of size $\le
\lambda$, then the second player has a winning strategy in the Ehrenfeucht-Fra\"\i ss\'e game of
length $\lambda^+$ on $\prod_i M_i/D$ and $\prod_i N_i/D$. If in addition $2^{\lambda}=\lambda^+$
and $i<\lambda$ implies $|M_i|+|N_i|\leq \lambda^+$ this means that the ultrapowers are
isomorphic. This settles negatively conjecture 18 in \cite{CK}.

\end{abstract}

The problem of when two elementarily equivalent structures have isomorphic ultrapowers was
prominent in the model theory of the 1960's. Keisler \cite{k} proved, assuming GCH, that
elementarily equivalent structures  have isomorphic ultrapowers. Keisler's proof depended on GCH
both on the question of existence of good ultrafilters and on limiting the size of the
ultraproducts. More exactly, Keisler considered a language of size $\lambda$, models $M$ of size
$\leq \lambda^+$ and a $\lambda^+$-good countably incomplete ultrafilter $D$ on $\lambda$. He
proved that $M^{\lambda}/D$ is $\lambda^+$-saturated. Under the weaker assumption that $D$ is
regular he proved that $M^{\lambda}/D$ is $\lambda^+$-universal, i.e. every $N\equiv
M^{\lambda}/D$ can be elementarily embedded into it.

Shelah \cite{MR45:6608} improved the result by eliminating GCH: two structures $M$ and $N$ are
elementarily equivalent if and only if, for some $\lambda$ and some regular ultrafilter $D$ on
$\lambda$ the structures $M^\lambda/D$ and $N^\lambda/D$ are isomorphic. This left open the
following question, asked by Chang and Keisler as Conjecture 18 in \cite{CK}:

\begin{quote} Let $M$ and $N$ be structures of cardinality $\leq \lambda$ in a language of size
$\leq \lambda$ and let $D$ be a regular ultrafilter over $\lambda$. If $M\equiv N$, then
$M^\lambda/D \cong N^\lambda/D$.
\end{quote}

The Conjecture is a natural one as most of the model theory regarding ultrapowers is centered on
the regular ultrafilters. Also at the time of Keisler's question GCH was generally considered a
reasonable assumption for the question.

Also the Conjecture is formulated for models of size $\leq \lambda$, perhaps for accidental
reasons, but it seems more natural if $M$ and $N$ have cardinality $\leq \lambda^+$.

Conjecture 19 of \cite{CK}, which we also address in this paper, is:

\begin{quote} If $D$ is a regular ultrafilter over $\lambda$, then for all infinite $M$,
$M^\lambda/D$ is $\lambda^{++}$-universal.
\end{quote}

In \cite{769} the authors proved that the transfer principle
$(\aleph_1,\aleph_0) \rightarrow (\lambda^+,\lambda)$ implies for
all regular filters $D$ on $\lambda$ \begin{itemize}
  \item[$(1)_D$] For all $M$ in a language of size $\leq \lambda$, $M^{\lambda}/D$
  is $\lambda^{++}$-universal.
  \item[$(2)_D$] If $M_i$ and $N_i$ are elementarily equivalent models of
               a language of size $\le \lambda$, then the second player
               has a winning strategy in the Ehrenfeucht-Fra\"\i ss\'e
               game of length $\lambda^+$ on $\prod_i M_i/D$ and
               $\prod_i N_i/D$.
\end{itemize}
Assuming $2^{\lambda} = \lambda^+$, $(2)_D$ is equivalent to:
 \begin{itemize}
  \item[$(2')_D$] For $M_i,N_i$ as in $(2)_D$ of cardinality
                  $\le \lambda^+$, $\prod_i M_i/D \cong \prod_i N_i/D$.
  \end{itemize}

We note that regularity is necessary for $(1_D)$. I.e.

\begin{lemma}

For any filter $D$ on an infinite cardinal $\lambda$, if for all infinite $N$ the structure
$N^{\lambda}/D$ is $\lambda^{+}$-universal, then $D$ is regular.

\end{lemma}

\begin{proof} For $k=1,2$ let ${\cal M}_k=\langle M_k, P^k_i\rangle_{i<\lambda}$, where
$M_1=\lambda$, $M_2=\lambda+1$ and  the $P^k_i$ are defined as
follows. Let $k=1$ and let $D_0$ be a regular filter on $\lambda$.
Let $\{A_i\}\subseteq D_0$ witness the regularity of $D_0$. Thus
for $\alpha<\lambda, \{i<\lambda|\alpha\in A_i\}$ is finite. Set
$P^1_i=A_i$. Now let $k=2$.  Let $P^2_i=A_i\cup\{\lambda\}$. Now
let $D$ be any ultrafilter on $\lambda$ and suppose there is an
elementary embedding $g$ mapping $M_2$ into $M_1^{\lambda}/D$. Let
$g(\lambda)=[f]$ and let $X_{\alpha}=\{i|f(i)\in P^1_{\alpha}\}$.
For each $\alpha<\lambda, a\in P^2_{\alpha}$ implies
$X_{\alpha}\in D.$ It is easy to see  that $\{X_{\alpha}\}_{\alpha
< \lambda}$ is a regular family in $D$. \qed

\end{proof}

We note that a similar argument can be used to show that $(1)_D$
fails if the language of $M$ has size $\lambda^+$.

This was a partial answer to the above Conjectures 18 and 19. In this paper we show the converse
for singular strong limit $\lambda$. Under GCH this is necessary as by Chang's Two Cardinal
Theorem the transfer principle $(\aleph_1,\aleph_0) \rightarrow (\lambda^+,\lambda)$ can only
fail, in the presence of GCH, for singular $\lambda$. It is known that $(\aleph_1,\aleph_0)
\nrightarrow (\aleph_{\omega+1},\aleph_{\omega})$ + GCH is consistent relative to the consistency
of a supercompact cardinal. It follows that the statement $(2)_D$is independent of ZFC relative to
the consistency of supercompact cardinals. On the other hand $(\aleph_1,\aleph_0) \rightarrow
(\lambda^+,\lambda)$ holds for all $\lambda$ if $V = L$.

In fact we show more. The results of \cite{769} were obtained
using a finitary version, denoted here by $\square^{fin}_{\lambda,
D}$, of the combinatorial principle $\square^{b^{*}}_{\lambda}$
from \cite{269}, which is equivalent for all $\lambda$ to
$(\aleph_1,\aleph_0) \rightarrow (\lambda^+,\lambda)$. We showed
in \cite{769} that whenever $\lambda$ is singular strong limit,
$\square^{fin}_{\lambda, D}$ is actually equivalent to
$\square^{b^{*}}_{\lambda}$, using Theorem 2.3 and Remark 2.5 of
\cite{269}. Thus $\square^{fin}_{\lambda, D}$ is equivalent to
$(\aleph_1,\aleph_0) \rightarrow (\lambda^+,\lambda)$, again for
$\lambda$ singular strong limit. The consistency of e.g.
GCH+$\neg\square^{b^{*}}_{\lambda}$ follows. Precisely we showed:

\begin{lemma} \cite{769} Let $D$ be a regular filter on $\lambda$ where
$\lambda$ is a singular strong limit cardinal. Then $\square^{fin}_{\lambda, D}$ is equivalent to
$\square^{b^{*}}_{\lambda}$.

\end{lemma}

The final piece needed for obtaining the transfer principle from
$(1)_D$, $(2)_D$ and $(2')_D$ for singular strong limit $\lambda$
involves proving their equivalence with the principle
$\square^{fin}_{\lambda, D}$, which equivalence is proved for all
$\lambda$ (Theorem~\ref{three}).

We note that $(2)_D$ is more robust than what was originally
conjectured, i.e. we have given a condition on when player II has
a winning strategy in the EF game of length $\lambda^+$ on the two
structures $\prod_i M_i/D$ and $\prod_i N_i/D$. Thus again if
$\lambda$ is singular strong limit, then the model theoretic
$(1)_D$ and $(2)_D$ are equivalent to the set theoretic
$(\aleph_1,\aleph_0) \rightarrow (\lambda^+,\lambda)$.

We need the following definition, from Lemma 4 of \cite{769}:

\def\uai{u^{\alpha}_i}
\def\uzi{u^{\zeta}_i}
\def\uei{u^{\epsilon}_i}
\def\ugi{u^{\gamma}_i}

\begin{definition}

Let $D$ be a regular filter on $\l$. If there exist sets \(\uzi\)
and integers $n_i$ for each $\zeta<\lambda^+$ and $i<\lambda$ such
that for each \(i,\zeta\)

\begin{description} \item[(i)] \(|\uzi|<n_i\) \item[(ii)]
\(\uzi\subseteq\zeta\) \item[(iii)] Let $B$ be a finite set of
ordinals and let $\zeta$ be such that $B\subseteq\zeta
<\lambda^+$. Then $\{i:B\subseteq \uzi\}\in D$ \item[(iv)]
Coherence: \(\gamma\in\uzi\Rightarrow\ugi=\uzi\cap\gamma\),

\end{description}

then we say that $\square^{fin}_{\lambda, D}$ holds.

\end{definition}

Our main result:

\begin{theorem}\label{three}
Assume $\lambda \ge \aleph_0$ and $D$ is a regular filter on
$\lambda$. Then the following conditions are equivalent:
\begin{itemize}
  \item[$(i)$] $\square^{fin}_{\lambda, D}$.
  \item[$(ii)$] If $M_i$ and $N_i$, $i < \lambda$, are elementarily equivalent models of
                a language of cardinality $\le \lambda$, then the second
                player has a winning strategy in the
                Ehrenfeucht-Fra\"{i}ss\'e game of length $\lambda^+$ on
                $\prod_i M_i/D$ and $\prod_i N_i/D$.
  \item[$(iii)$] If $M$ and $N$ are structures of a language of
                 cardinality $\le \lambda$, $N \equiv M^{\lambda}/D$
                 and $|N|\leq\lambda^+$,
                 then there is a homomorphism $N \rightarrow
                 M^{\lambda}/D$.
  \item[$(iv)$] If $\Delta$ is a set of quantifier-free formulas and
                $M^{\lambda}/D$ satisfies every existential
                $\Delta$-sentence (i.e. a sentence of the form $\exists
                \vec{x}(\phi_1\wedge...\wedge\phi_n)$, where each $\phi_i$ is in
                $\Delta$) true in $N$, $|N|\leq \lambda^+$, then there is a
                $\Delta$-homomorphism $N \rightarrow M^{\lambda}/D$, i.e.
                a homomorphism $N \rightarrow M^{\lambda}/D$ which preserves
                $\Delta$ formulas. \end{itemize}
                Additionally, if $D$ is an ultrafilter, then $(i)-(iv)$ are equivalent to

                \begin{itemize}
  \item[$(v)$] If $M$ is a structure in a language of cardinality
               $\le \lambda$, then $M^{\lambda}/D$ is
               $\lambda^{++}$-universal.
\end{itemize}

Moreover in $(ii)$, $(iii)$ and $(iv)$ we can equivalently assume that the models $M_i$, $N_i$,
$M$ and $N$ have cardinality $\le \lambda^+$.

\end{theorem}

\begin{proof}
$(i) \rightarrow (ii)$, $(i) \rightarrow (iii)$ and $(i)
\rightarrow (iv)$ follow from the "$\Delta$-existential" version
of Theorem~2 of \cite{769} which gives a $\Delta$-homomorphism of
$N$ into $M^{\lambda}/D$ for any set $\Delta$ of first-order
formulas such that every $\Delta$-existential sentence true in $N$
is true in $M$. If $D$ is an ultrafilter, \cite{769} gives $(i)
\rightarrow (v)$ and, on the other hand, $(v) \rightarrow (iii)$
is straightforward.

$(ii) \rightarrow (iii)$: It follows from $N \equiv M^{\lambda}/D$
that there are $M_0 \equiv M$ and a homomorphism $N \rightarrow M_0$
(i.e. a mapping from $N$ to $M_0$ which respects the functions and
relations of $N$) such that $M_0$ has cardinality $\le \lambda^+$. By
$(ii)$ the second player has a winning strategy in the
Ehrenfeucht-Fra\"{i}ss\'e game of length $\lambda^+$ on
$M^{\lambda}_{0}/D$ and $M^{\lambda}/D$. Using this winning strategy
we get easily a homomorphism $N \rightarrow M^{\lambda}/D$, by
composing the appropriate mappings.

$(iii) \rightarrow (iv)$: Let $M$ and $N$ be as in $(iv)$. Let
$N^*$ be an expansion of $N$ obtained by giving a name to every
$\Delta$-definable relation. Let $(M^{\lambda}/D)^*$ be obtained
similarly from $M^{\lambda}/D$. Let $M^{*}_{0} \equiv
(M^{\lambda}/D)^*$ such that there is a homomorphism $N^*
\rightarrow M^{*}_{0}$. (The existence of such $M^*_0$ follows
from the fact that we can find $M^{*}_{0} \equiv
(M^*)^{\lambda}/D$ and a homomorphism $N^* \rightarrow M^{*}_{0}$.
But then $(M^*)^{\lambda}/D$ is canonically embeddable into
$(M^{\lambda}/D)^*$). By $(iii)$ there is a homomorphism
$M^{*}_{0} \rightarrow M^{\lambda}/D$. Thus there is a
$\Delta$-homomorphism $N \rightarrow M^{\lambda}/D$.

$(iv) \rightarrow (i)$:  Let $\lambda, D$ be given and let $\tau$
denote a language of cardinality $\leq \lambda$. It suffices to
prove the following
\smallskip

\noindent{\bf Claim.} There exist $M, N$ such that \begin{itemize}
  \item[$a)$] $\lvert M \rvert = \lambda, \lvert N \rvert = \lambda^+$
  \item[$b)$] $\tau_M = \tau_N$ and $\lvert \tau_N \rvert \le \lambda$
  \item[$c)$] $M \equiv N$
  \item[$d)$] For $\Delta =$ the quantifier free formulas of $\tau_M$,
  $N$ has a $\Delta$-homomorphism into
              $M^{\lambda}/D$, and hence $\square^{fin}_{\lambda, D}$ holds.
               \end{itemize}

 % \end{claim}

\begin{proof}
Let $\tau^* = \{ F_{\alpha} \mid \alpha < \lambda \} \cup \{ <
\}$, for $F_{\alpha}$ a unary function symbol. Let $K$ be the
family of all structures $M$ such that \begin{itemize}
  \item[$K1)$] $M$ is a finite $\tau^*$-structure.
  \item[$K2)$] The universe of $M$ is $\{ 0,1, \ldots, k-1\}$,
  for some $k\in\mathbb{N}$, $k\neq 0$, denoted $\eta(M)$.
  \item[$K3)$] $M \models \forall x (F_{\alpha}(x) \le x)$ for all
              $\alpha < \lambda$.
  \item[$K4)$] If $m_1 = F^{M}_{\alpha_1}(m)$, $m_2 = F^{M}_{\alpha_2}(m)$
              and $m_1 < m_2$, then there exists a $\beta < \lambda$ such that
              $m_1 = F^{M}_{\beta}(m_2)$.
  \item[$K5)$] If $F^{M}_{\alpha_2}(m_3) = m_2$, $F^{M}_{\alpha_1}(m_2) = m_1$
              and $m_1 < m_2 < m_3$, then there exists $\alpha_3 < \lambda$
              such that $F^{M}_{\alpha_3}(m_3) = m_1$.
  \item[$K6)$] $w(M) =_{df} \{ \alpha \mid F^{M}_{\alpha} \text{ is not the
              identity} \}$ is finite.

  \item[$K7)$] For $m_1 < m_2 < k$ there is exactly one $\alpha$ such that $m_1 =
  F_{\alpha}(m_2))$.

\end{itemize}
We note that $K$ is non-empty, taking $K$ to be, e.g., a one element
structure. Let $\{ M_i \mid i < \lambda\}$ list $K$. We will add the
$M_i$ together into one structure. I.e., we define a model $M^*$ for
$\tau = \tau^* \cup \{E\}$ such that

\begin{itemize}
  \item[$a_K)$] The universe of $M^* = \cup\{ \{ i\} \times M_i \mid i < \lambda \}$
  \item[$b_K)$] $E^{M^*} = \{ \langle (i_1,m_1),(i_2,m_2)\rangle \mid
              m_1 < \eta(M_{i_1}), m_2 < \eta(M_{i_2}) \text{ and }
              i_1=i_2 \}$
  \item[$c_K)$] $<^{M^*} = \{ \langle (i,m_1),(i,m_2) \rangle \mid
              m_1 < m_2 < \eta(M_{i}) \}$
  \item[$d_K)$] $F^{M^*}_{\alpha}(i,m) = \langle i,F^{M_i}_{\alpha}(m) \rangle$.
  \end{itemize}

Now for $\rho < \lambda^+$ let $h_{\rho}$ be a partial one to one function from $\lambda$ onto
$\rho$, and let $\langle a_{\rho} \mid \rho < \lambda^+ \rangle$ be a set of new constant symbols.
\bigskip

\noindent{\bf Subclaim.\ } There is $N^*$ such that
\begin{itemize}
  \item[$a)$] $N^*$ is a $\tau^{**}$ structure of cardinality $\lambda^+$, where $\tau^{**} = \tau \cup
              \{ a_{\rho} \mid \rho < \lambda^+ \} \cup
              \{ \overline{m} \}_{m \in M}$,
  \item[$b)$] $M^* \preccurlyeq N^* \restriction
   \tau$.
  \item[$c)$] $a_{\rho} E^{N^*} a_0$ for $\rho<\lambda^+$,
  \item[$d)$] $a_{\rho} <^{N^*} a_{\xi}$, for $\rho < \xi <
  \lambda^+$,
  \item[$e)$] $N^* \models F_j(\bar{a}_{\rho}) =
  \bar{a}_{\epsilon}$, if $h_{\rho}(j) =
  \epsilon$, for $j<\lambda$ and $\epsilon < \rho < \lambda^+$. \end{itemize}
\bigskip

\noindent{\bf Proof of Subclaim.} Let $T = Th(M^*,\overline{m})_{m\in
M^*} \cup \{ a_{\rho} E a_0 \}_{\rho<\lambda^+} \cup \{ a_{\rho} <
a_{\xi} \mid \rho < \xi < \lambda^+ \}\cup \{
F_j(a_{\rho})=a_{\epsilon}\mid h_{\rho}(j)=\epsilon, \epsilon < \rho
< \lambda^+\}$. We claim that $T$ is consistent. To see this, let
\[T_0= \{ \varphi_i(\vec{\overline m}\}_{i=1,\ldots, n} \cup
\{a_{\rho_i}Ea_0\}_{i=1,\ldots k} \cup \{a_{\epsilon_i}<a_{\zeta_i}
\}_{i=1,\ldots,l} \cup
\{F_{\alpha_i}(a_{\mu_i})=a_{\nu_i}\}_{i=1,\ldots,m}\] be a finite
part of $T$, where $a_{i_1}, \ldots, a_{i_k}$ and
$\overline{m}_{j_1}, \ldots, \overline{m}_{j_l}$ are all the
parameters occurring in $T_0$. Let $I_0=\{i\mid $ for some
$j=j_1\ldots, j_l$, $m_j = (i,a)$, $i < \lambda$, $a < \eta(M_i)\}$.
We can find $M_{i_0}$, $i_0\in \lambda \smallsetminus I_0$, such that
if we interpret the constants $a_{\rho}$ in $M_{i_0}$ and $\overline
m_{j_i}$ by $m_{j_i}\in M^*$, then this expansion of $M^*$ is a model
of $T_0$ and thus $T_0$ is consistent. The subclaim is proved.

%To this end, define the domain of $N_0$ to be the ordered segment \{0,1,\ldots,k-1\}. Define
%$a^{N_0}_{i_r}=r-1, r=1,\ldots,k$ and define $F^{N_0}_j(a^{N_0}_{\rho})=a^{N_0}_{\epsilon}$ if
%$j,\rho,\epsilon$ respectively $=\alpha_i,\mu_i, \nu_i$, for $i=1,\ldots, m$. Note that each
%$F_j^{N_0}$ so defined is a regressive function. Inductively define new functions
%$F^{N_0}_{\beta}$ as follows: if $a^{N_0}_{i_r}$ is not already related to $a^{N_0}_{i_s}$ by some
%$F^{N_0}_j$ already defined, then for each such $a_{i_r}, r=2,\ldots,k$, define
%$F_{\beta_s}(a_{i_r})=a_{i_s}$, where $s=1,\ldots,r-1$ and $\beta_s$ is a new ordinal not
%contained in $T_0$ and not equal to any previously defined $\beta_{s'}$. (I.e. we "connect" each
%element $a$ to all predecessors it is not already connected to via a new function symbol.)
%Finally, for all "remaining" $\alpha$, i.e. all $\alpha$ such that $\alpha$ does not occur in
%$T_0$ and $\alpha\neq\beta_s$, for any newly defined $F_{\beta_s}$, define $F_{\alpha}$ to be the
%identity. We claim that $N_0\in K$, and moreover that $N_0\neq M_i$, for any $i\in I_0$. $K1$ and
%$K2$ are obvious. $K3$ holds as we have already noted that all the defined non-identity functions
%are regressive. $K4$, $K5$ and $K6 $ follow from the construction of $N_0$. $K7$ follows from the
%fact that the $H_{\rho}$'s are one to one. The subclaim is proved.

\medskip

Now let $N^*$ be as in claim 1 and let $N=N^*\restriction \tau_M^*$.
We note that the pair of structures $N$ together with the $M$ defined
above satisfying $a_K)-d_K)$, satisfy the hypothesis of Theorem
2-(iv), i.e. $(M^*)^{\lambda}/D$ satisfies every existential
$\Delta$-sentence true in $N$ where $\Delta =$ the quantifier-free
formulas of $\tau_M^*$. This is because $N \equiv M^*$ and these
$\Delta$-sentences are preserved under reduced products. Therefore by
$(iv)$ there is a $\Delta$-homomorphism $g \colon N \rightarrow
(M^*)^{\lambda}/D$. Let $g(a_{\rho})$ be denoted by $f_{\rho}/D$.

\medskip

We are now ready to define the sets $\uzi$ referred to in (i)-(iv) of
the condition $\square^{fin}_{\lambda,D}$. To this end, for $\epsilon
< \rho < \lambda^+$, $\rho \ge \lambda$, define \[ A_{\epsilon,\rho}
= \{ j < \lambda \mid M^* \models (f_{\rho}(j) >
f_{\epsilon}(j))\wedge F_{i(\epsilon,\rho)}(f_{\rho}(j)) =
f_{\epsilon}(j) \}. \] Note that if $\epsilon < \rho < \lambda^+$,
$\rho \ge \lambda$, then $A_{\epsilon, \rho} \in D$, since
$(M^*)^{\lambda}/D\models (g(a_{\epsilon})<g(a_{\rho}))\wedge
F_{i(\epsilon, \rho)}(g(a_{\rho}))=g(a_{\epsilon})$.
\medskip
For each $\rho \in (\lambda,\lambda^+)$ and $j \le \lambda$ define
\[ W_{\rho,j} = \{ \epsilon < \rho \mid \epsilon \ge \lambda
\text{ and } j \in A_{\epsilon,\rho} \}. \]

%We claim that if $u^{\rho}_j =_{df} W_{\rho,j}$ then
%$\square^{fin}_{\lambda,D}$ holds.

\medskip

%Proof of Claim. First note that under the $\Delta$-homomorphism $g
%\colon N \rightarrow M^{\lambda}/D$, if $M^j$ ($= M$) is the
%$j$-th "component" of the reduced product $M^{\lambda}/D$, then
%there is a unique $i_0=i_0(j)$, $i_0<\lambda$, such that each
%image element $f_{\rho}(j)/D$ belongs to the same equivalence
%class $i_0 \times M_{i_0}$ of $M$. This justifies introducing the
%following notation for the proofs of (i)-(ii) and (iv) of
%$\square^{fin}_{\lambda, D}$: for fixed $j$ and $i_0=i_0(j)$ as
%above, if $f_{\rho}(j) = (i_0,m(i_0,\rho))$, denote $m(i_0, \rho)$
%by $\tilde{f_{\rho}(j)}$.

\noindent First note that without loss of generality
\smallskip

$(*)_1$:  $M^*\models f_{\rho}(j) E f_0(j)$ for every $j,\rho$.

\smallskip

\noindent We claim that

\smallskip

$(*)_2$:  if we choose $u^{\rho}_j =_{df} W_{\rho,j}$ then $\langle
u^{\rho}_j | \rho < \lambda^+,
j < \lambda \rangle $ exemplifies $\square^{fin}_{\lambda,D}$, i.e.
these objects satisfy the
demands $(i)-(iv)$ of Definition 3.
\smallskip

\noindent By $(*)_1$ we can let $f_{\rho}(j)= (i_j, m(i_j,\rho)).$
Clearly if we prove $(*)_2$ we
are done.

\smallskip

(i), (ii): $W_{\rho,j}$ is a finite subset of $\rho$:
\begin{eqnarray*}  \{ \epsilon < \rho \mid \epsilon \ge \lambda
\text{ and } j \in A_{\epsilon,\rho} \}  &= \{ \epsilon < \rho
\mid \epsilon \ge \lambda \text{ and } M \models
    F_{i(\epsilon,\rho)}(f_{\rho}(j)) = f_{\epsilon}(j) \} \\  &= \{ \epsilon < \rho \mid \epsilon \ge
    \lambda \text{ and } M_{i_0} \models
    F^{M_{i_0}}_{i(\epsilon,\rho)}(f_{\rho}(j)) = f_{\epsilon}(j)\}. \end{eqnarray*} But $w(M_{i_0})$ is finite, and therefore so is $W_{\rho,j}$. Thus if $n_i$ is taken to be $w(M_{i_o})$ then (i) and (ii) of $\square^{fin}_{\lambda,D}$ are satisfied.

(iv): (coherency) holds of $W_{\rho,j}$, i.e. if $\rho_1 < \rho_2 <
\lambda^+$ and $j < \lambda$ is given, if $\rho_1 \in W_{\rho_2,j}$
then $W_{\rho_1,j} = W_{\rho_2,j} \cap \rho_1$. Why? Let $\alpha \in
W_{\rho_1,j}$. Then $\lambda\leq\alpha < \rho_2$ and $j \in
A_{\alpha,\rho_1}$. But then $M^*\models
F_{i(\alpha,\rho_1)}(f_{\rho_1}(j)) = f_{\alpha}(j)$ and therefore
$M_{i_0}\models
F^{M_{i_0}}_{i(\alpha,\rho_1)}(f_{\rho_1}(j))=f_{\alpha}(j)$. $\rho_1
\in W_{\rho_2,j}$ and therefore $j \in A_{\rho_1,\rho_2}$. This means
$M_{i_0} \models F^{M_{i_0}}_{i(\rho_1,\rho_2)}(f_{\rho_2}(j)) =
f_{\rho_1}(j)$. By the definition of $M_{i_0}$, there is $\beta$ such
that $M_{i_0} \models F^{M_{i_0}}_{\beta}(f_{\rho_2}(j)) =
f_{\alpha}(j)$. But $\beta = i(\alpha,\rho_2)$ and $j \in
A_{\alpha,\rho_2}$ and therefore $\alpha \in W_{\rho_2,j}$. For the
other direction suppose $\alpha \in W_{\rho_2,j} \cap \rho_1$. Then
$\alpha \ge \lambda$ and $j \in A_{\alpha,\rho_2}$, i.e.\ $M_{i_0}
\models F^{M_{i_0}}_{i(\alpha,\rho_2)}(f_{\rho_2}(j)) =
f_{\alpha}(j)$. $\rho_1 \in W_{\rho_2,j}$ means that as before $j \in
A_{\rho_1,\rho_2}$, i.e.\ $M_{i_0} \models
F^{M_{i_0}}_{i(\rho_1,\rho_2)}(f_{\rho_2}(j)) = f_{\rho_1}(j)$. But
then since $\alpha < \rho_1$, there is $\beta$ such that $M_{i_0}
\models F^{M_{i_0}}_{\beta}(f_{\rho_1}(j)) = f_{\alpha}(j)$, i.e.\ as
before $j \in A_{\rho_1,j}$ and $\alpha \in W_{\rho_1,j}$.

To see that (iii) is satisfied, let $B\subseteq \lambda^+$ be a
finite set of ordinals such that $B\cap\lambda=\emptyset$, and let
$\rho$ be such that $B \subseteq \rho < \lambda^+$. We wish to show
that $\{ j \mid B \subseteq W_{\rho,j} \} \in D$. Let
$B=\{j_1,\ldots,j_n\}$. Recall that by (e) of the subclaim,
$N^*\models F_{i(j,\rho)}(a_{\rho})=a_j$ whence \[C_k =\{i \mid
M^*\models F_{i(j_k,\rho)}(f_{\rho}(i))= f_{j_k}(i)\}\in D,\] for
each $k=1,\ldots,n$. Also by definition if $i\in C_k$ then $i\in
A_{j_k,\rho}$ whence $j_k\in W_{\rho,i}$. Thus $C_1\cap \cdots \cap
C_n \subseteq \{ i \mid B \subseteq W_{\rho,j} \}\in D.$ The claim is
proved.

Now if we transfer the sets $u^{\rho}_j$, $\lambda \leq \rho <
\lambda^+$ to all of $\lambda^+$, (iv) implies (i), proving the
Claim.

The claim concerning the equivalent formulations involving models of size $\le \lambda^+$ follows
from the fact that in the derivation $(iv) \rightarrow (i)$ above we have $|M|=\lambda$ and
$|N|=\lambda^+$. Thus Theorem 4 is proved. \qed

\end{proof} \end{proof}

\begin{corollary}
Assume $\lambda \geq \aleph_0$, $2^{\lambda}=\lambda^+$ and $D$ is a regular filter on $\lambda$.
Then the following are equivalent:

\begin{itemize}
  \item[$(i)$] $\square^{fin}_{\lambda, D}$.
  \item[$(ii)$] If $M_i$ and $N_i$, $i < \lambda$, are elementarily equivalent models of
                a language of cardinality $\le \lambda$, and $|M_i|$, $|N_i| \le \lambda^+$ then
                $\prod_i M_i/D \cong \prod_i N_i/D$.
\end{itemize}

\end{corollary}

\begin{corollary}
 GCH $+$ the failure of properties (i)-(iv) of Theorem 4 for
$\lambda=\aleph_{\omega}$ is consistent relative to the
consistency of supercompact cardinals.

\end{corollary}

\proof Assume GCH and $(\aleph_1,\aleph_0) \nrightarrow
(\aleph_{\omega+1},\aleph_\omega)$. It is well-known (cf.
\cite{CK} Proposition 4.3.5) that there is a regular ultrafilter
$D$ on \(\aleph_\omega\). The principle
$\square^{b^{*}}_{\aleph_\omega}$ fails. Thus
$\square^{fin}_{\aleph_\omega,D}$ fails by \cite{769}. By
Theorem~\ref{three}, (i)-(v) fail for the regular ultrafilter $D$.
\qed

A drawback of Corollary 6 is that it deals with filters rather than ultrafilters, which was
originally the most interesting case. This case will be dealt with in a work in preparation.

%The ultrafilter case of the above Corollary is the subject of
%ongoing work.

 \end{document}